\documentclass{amsart}
\usepackage[utf8]{inputenc}
\usepackage{enumerate}
\usepackage{graphicx}
\usepackage{wrapfig}
\usepackage{subcaption}
\usepackage{listings}
\usepackage{color}
\definecolor{mygreen}{RGB}{28,112,30} % color values Red, Green, Blue
\definecolor{mylilas}{RGB}{170,55,241}
\definecolor{myblue}{RGB}{20,30,171}
\definecolor{myred}{RGB}{200,10,30}
\usepackage{microtype}
\usepackage{multirow}
\usepackage{array}
\newcounter{rowcount}
\usepackage{amsmath}
\usepackage{amsfonts}
\usepackage{amssymb}
\usepackage{amscd}
\usepackage{amstext}
\usepackage{amsbsy}
\usepackage{amsopn}
\usepackage{amsthm}
\usepackage{thmtools}
\usepackage{amsxtra}

\usepackage{float}
\usepackage{bbm}
\usepackage{url}
\usepackage{hyperref}% wants to be last
\hypersetup{
colorlinks=true,
linkcolor=myred,
citecolor=blue}
\newtheoremstyle{mytheoremstyle} % name
  {12pt}                    % Space above
  {\topsep}                    % Space below
  {\slshape}                   % Body font
  {}                           % Indent amount
  {\bfseries}                   % Theorem head font
  {.}                          % Punctuation after theorem head
  {.5em}                       % Space after theorem head
  {}  % Theorem head spec (can be left empty, meaning ‘normal’)

\newtheoremstyle{mydefstyle} % name
  {12pt}                    % Space above
  {12pt}                    % Space below
  {}                   % Body font
  {}                           % Indent amount
  {\bfseries}                   % Theorem head font
  {.}                          % Punctuation after theorem head
  {.5em}                       % Space after theorem head
  {}  % Theorem head spec (can be left empty, meaning ‘normal’)

  \newtheoremstyle{plainsl}%
  	{}
  	{\topsep}
  	{\slshape} % only non-default setting
  	{}
  	{\normalfont\bfseries}
  	{.}
  	{ }
  	{}

\theoremstyle{plainsl}
\newtheorem{thm}{Theorem}[section]
\newtheorem*{thm*}{Theorem}

\newtheorem{cor}[thm]{Corollary}
\newtheorem{lem}[thm]{Lemma}
\newtheorem*{lem*}{Lemma}

\theoremstyle{remark}

\theoremstyle{mydefstyle}

\newtheorem{xmpl}[thm]{Example}
\newtheorem*{note*}{Remark}

\newcommand{\Z}{\mathbb{Z}}

\newcommand{\R}{\mathbb{R}}

\newcommand{\Cx}{\mathbb{C}}

\newcommand{\C}{\mathcal{C}}

\newcommand{\eps}{\varepsilon}

\renewcommand{\phi}{\varphi}
\newcommand{\sbeq}{\subseteq}

\newcommand\comp[1]{{\mkern2mu\overline{\mkern-2mu#1}}}
\newcommand\diff{\mathbin{\mkern-1.5mu\setminus\mkern-1.5mu}}% for \setminus

\renewcommand\qed{%
	\ifmmode\eqno\sqr53
	\else\nolinebreak\ \hfill\sqr53\medbreak\fi}
\renewcommand\proof{\noindent\textsl{Proof. }}
\newcommand\sqr[2]{{\vbox{\hrule height.#2pt
    \hbox{\vrule width.#2pt height#1pt \kern#1pt
        \vrule width.#2pt}\hrule height.#2pt}}}

\DeclareMathOperator\Aut{Aut}

\DeclareMathOperator\tr{Tr}

\title{Strongly cospectral vertices in normal Cayley graphs}
\author{Arnbjörg Soffía Árnadóttir \& Chris Godsil}
\thanks{Both authors acknowledge the support of C.\ Godsil's NSERC (Canada), Grant No.\ RGPIN-9439}
\date{\today}

\begin{document}
\renewcommand{\itshape}{\slshape}

\begin{abstract}
  We prove an upper bound on the number of pairwise strongly cospectral vertices in a normal Cayley graph, in terms of the multiplicities of its eigenvalues. We use this to determine an explicit bound in Cayley graphs of $\Z_2^d$ and $\Z_4^d$. We also provide some infinite families of Cayley graphs of $\Z_2^d$ with a set of four pairwise strongly cospectral vertices and show that such graphs exist in every dimension.
\end{abstract}

\maketitle
\setlength{\parskip}{0cm}
\tableofcontents
\setlength{\parskip}{0.2cm}
\section{Introduction}

The spectrum of a graph (meaning the spectrum of its adjacency matrix) carries a lot of information about the graph itself. It is therefore natural to ask: when do two graphs have the same spectrum? Such graphs are said to be \textsl{cospectral}.

A related concept is cospectrality of vertices. Two vertices, $u$ and $v$ in a graph $X$ are \textsl{cospectral} if the graphs $X\diff u$ and $X\diff v$ are cospectral.
In this paper, we focus on a stronger property of vertices in a graph called strong cospectrality. Let $X$ be a graph with adjacency matrix $A$. For each eigenvalue, $\theta_r$ of $A$, let $E_r$ denote the projection onto the eigenspace corresponding to $\theta_r$. Further, denote by $e_u$ the standard basis vector indexed by the vertex $u$. We say that vertices $u$ and $v$ are \textsl{parallel} if the vectors $E_re_u$ and $E_re_v$ are parallel for all $r$. If $u$ and $v$ are both cospectral and parallel, we call them \textsl{strongly cospectral}.

Strong cospectrality of vertices has some combinatorial implications for the graph. For instance, if $u$ and $v$ are strongly cospectral in $X$, then an automorphism of $X$ that fixes $u$ also fixes $v$. Further, as we show in Section \ref{sec:blocks} of this paper, a maximal set of pairwise strongly cospectral vertices forms a block of imprimitivity under any transitive automorphism group of $X$. In unpublished lecture notes \cite[Theorem 6.8.2]{Coutinho2021}, Coutinho and Godsil characterize strongly cospectral vertices, providing five equivalent definitions. This includes connections to walk matrices and $\R[A]$-modules.
The initial motivation, however, for studying strong cospectrality comes from physics, in the form of so-called state transfer.

Let $X$ be a graph with adjacency matrix $A$. The \textsl{continuous-time quantum walk} on $X$ at time $t$ is given by the matrix $U(t) :=e^{itA}$. We call $U(t)$ the \textsl{transition matrix} of the quantum walk. Let $u$ and $v$ be vertices in $X$. We say that there is \textsl{perfect state transfer} from $u$ to $v$ at time $t$ if $|U(t)_{u,v}|=1$.
A quantum walk can be thought of as the journey of a quantum particle around a graph and perfect state transfer describes a scenario where the particle has completely moved from one quantum state to another. Therefore, perfect state transfer is of significant interest in quantum physics and quantum computing.

The concept of perfect state transfer was initiated by Bose in 2003 \cite{bose2003}. Since then, many others have studied this phenomena, including Christandl et al \cite{christandl2004,christandl2005} and Kay \cite{kay2010,kay2011}.
In \cite[Lemma 11.1]{godsil2012state}, Godsil showed that strong cospectrality between vertices is a necessary condition for perfect state transfer to occur between them.

Strong cospectrality, however, is a weaker property than perfect state transfer. Whereas perfect state transfer can only occur between a pair of vertices, a vertex can be strongly cospectral to more than one other vertex. The smallest example of this is the Cartesian product of $P_2$ and $P_3$ (the paths on two and three vertices respectively), in which the four vertices of degree two are pairwise strongly cospectral. We are interested in knowing how large such sets of vertices can be in certain types of graphs.

In their 2017 paper \cite[Lemma 10.1]{godsil2017strongly}, Godsil and Smith showed that if all vertices in a graph $X$ are pairwise strongly cospectral, then $X=K_2$. This is not true for cospectrality of vertices, since in any vertex-transitive graph, all the vertices are cospectral. In the same paper, they ask the question of whether there exists a tree with at least three pairwise strongly cospectral vertices. It has been shown very recently, by Coutinho, Juliano and Spier, that such trees do not exist \cite{Coutinho2022trees}.

In this paper, we consider sets of pairwise strongly cospectral vertices in Cayley graphs having the property that the connection set is closed under conjugation. We will call such Cayley graphs \textsl{normal}. Note that in the literature, there are two non-isomorphic definitions of a normal Cayley graph; this one is consistent with a paper by Larose et al from 1998 \cite{Larose1998}.

All Cayley graphs are vertex transitive, and so any two vertices in a Cayley graph are cospectral. They are therefore a somewhat natural choice when looking for strongly cospectral vertices. The motivation for studying normal Cayley graphs in particular comes from the fact that their spectrum can be calculated using the irreducible characters of the group, as we explore in Section \ref{sec:normalspec}. This property will prove extremely useful in our investigation. It turns out that despite the high symmetry of Cayley graphs we cannot have large sets of vertices that are pairwise strongly cospectral, at least in the normal case.

Our main results are the following.
We prove that if $X$ is a normal Cayley graph, and $m$ is the largest multiplicity of an eigenvalue of $X$, then the number of pairwise strongly cospectral vertices is bounded above by $|V(X)| / m$ (Theorem \ref{thm:multbound}). We apply this result to cubelike graphs (Cayley graphs of $\Z_2^d$), showing that such a graph has at most
\[2^{\lceil d/2\rceil-1} < \sqrt{2^d}\]
pairwise strongly cospectral vertices (Corollary \ref{cor:cubebound}). We also provide examples of this bound being tight for $d=3,4,5,6$. Theorem \ref{cor:Z4dbound} states a similar bound for Cayley graphs of $\Z_4^d$. Finally, we show in Theorems \ref{thm:exodd} and Theorem \ref{thm:exeven} that cubelike graphs with sets of four pairwise strongly cospectral vertices exist in every dimension $d\geq 5$.
We use this to prove Theorem \ref{thm:fourall} stating that there exist Cayley graphs on $n$ vertices with such sets of size four, for every integer $n$ divisible by $32$.

\section{Preliminaries}\label{sec:prelim}
Let $X$ be a simple, undirected graph and $A$ its adjacency matrix. Since $A$ is real and symmetric, it is diagonalizable and has real eigenvalues. We will refer to the eigenvalues and eigenvectors of the adjacency matrix as the eigenvalues and eigenvectors of the graph. Let $\theta_1,\dots,\theta_d$ be these (distinct) eigenvalues and let $E_r$ be the orthogonal projection onto the eigenspace of $\theta_r$. Then each $E_r$ is a polynomial in $A$ and we have $E_r^2=E_r$ for all $r$ and $E_rE_s=0$ if $r\neq s$. Moreover, $\sum_{r=0}^d E_r=I$ and we can write
\[A = \sum_{r=0}^d\theta_rE_r.\]
This is called the \textsl{spectral decomposition} of $A$ and the matrices $E_r$ are the \textsl{spectral idempotents} of $A$ (and of $X$).

Two graphs are said to be \textsl{cospectral} if they have the same spectrum. Vertices $u$ and $v$ in a graph $X$ are \textsl{cospectral} if the graphs $X\diff u$ and $X\diff v$ are cospectral. Let $e_u$ and $e_v$ denote the standard basis vectors indexed by $u$ and $v$, respectively.
Then, $u$ and $v$ are called \textsl{parallel} if for each $r=1,\dots,d$, the projections $E_re_u$ and $E_re_v$ are parallel. We say that $u$ and $v$ are \textsl{strongly cospectral} if they are both cospectral and parallel; equivalently, if for all $r$ we have $E_re_u = \pm E_re_v$. A set of pairwise strongly cospectral vertices that has size at least two will be called a \textsl{strongly cospectral set}.

Let $G$ be a group and $\C\sbeq G\diff 1$ an inverse-closed subset. The \textsl{Cayley graph} of $G$ with respect to the set $\C$ is the graph with vertex set $G$ where vertices $g$ and $h$ are adjacent if $hg^{-1}\in \C$. We denote this graph by $X(G,\C)$ and refer to $\C$ as its \textsl{connection set}. A Cayley graph is called \textsl{normal} if the connection set is a union of conjugacy classes of $G$, equivalently if $g^{-1}\C g = \C$ for all $g$ in $G$.
If $G$ is abelian, we refer to $X(G,\C)$ as a \textsl{translation graph}, and if $G$ is an elementary abelian $2$-group, we call the graph a \textsl{cubelike graph}.

\section{Group actions and blocks}\label{sec:blocks}

Let $G$ be a group acting transitively on a set $\Omega$ and denote the image of $\alpha\in\Omega$ under $g\in G$ by $\alpha^g$. A \textsl{block (of imprimitivity)} is a subset $B$ of $\Omega$ satisfying that for every $g\in G$, either $B^g=B$ or $B^g\cap B=\emptyset$. For further preliminaries in the theory of group actions, we refer the reader to \cite{dixon1996}.

We will show that strongly cospectral sets in a graph are blocks of imprimitivity under the action of any transitive automorphism group.

\begin{lem}\label{lem:blocks}
  Let $X$ be a graph and $G\leq\Aut(X)$ a group acting transitively on its vertices. Then, a maximal strongly cospectral set is a block under this action.
\end{lem}
\proof
  Let $A$ be the adjacency matrix of $X$ and let $E_0,\dots,E_d$ denote its spectral idempotents.
  Let $B\subseteq V(X)$ be a maximal strongly cospectral set and let $g\in G$. Since $g$ is a permutation of the vertices of $X$, we can think of it as a permutation matrix, $P_g$, mapping the standard basis vector $e_v$ to $e_{v^g}$ for all $v\in V(X)$.
  It is well known that a permutation matrix commutes with $A$ if and only if it is an automorphism of $X$. Further, since the spectral idempotents are polynomials in $A$, this implies that $P_g$ commutes with each of them.

  Suppose that there is some vertex $u$ in $B$ such that $u^g\in B$, and let $v\in B$. Then for all $r$ we have
  \[E_re_{v^g} = E_rP_ge_v = P_gE_re_v = \pm P_gE_re_u = \pm E_rP_ge_u = \pm E_re_{u^g}\]
  and so $v^g$ is strongly cospectral to $u^g$ and therefore $v^g\in B$. It follows that $B$ is a block of imprimitivity under the action of $G$.
\qed

Let $G$ be a group and $X=X(G,\C)$ a Cayley graph. For each $g\in G$, the map $V(X)\to V(X)$ defined by $x\mapsto xg$ is a graph automorphism of $X$, and so $G\leq \Aut(X)$. This action of $G$ on the vertices of the graph is regular, so in particular, Cayley graphs are vertex transitive.

\begin{lem}\label{lem:blockgroups}
  The blocks of the regular action mentioned above are precisely the subgroups of $G$ and their cosets.
\end{lem}
\proof
  If $H$ is a subgroup of $G$, the image of $H$ under $g\in G$ is the right coset $Hg$, which we know to be either equal to $H$ or disjoint from it. Conversely, if $B'\subseteq G$ is a block, since the action is transitive, some translate, $B$ of $B'$ will contain the identity. Let $x,y\in B$. Then $y=1^y\in B^y$, so $B^y=B$, but this implies that $x^y = xy\in B$. Similar argument shows that $x^{-1}\in B$ and so $B$ is a subgroup of $G$.
\qed

\begin{cor}\label{cor:blockgroups}
  Let $X=X(G,\C)$ be a Cayley graph. Under the action of $\Aut(X)$ on $X$, every block is a right coset of a subgroup of $G$. In particular, every block that contains the identity is a subgroup of $G$.
\end{cor}
\proof
  It is not too hard to see that if $H_1\leq H_2$ are groups acting on a set, then any block under the action of $H_2$ is also a block under the action of $H_1$. Therefore, since $G\leq \Aut(X)$, the result follows from Lemma \ref{lem:blockgroups}. \qed

By Lemma \ref{lem:blocks} and Corollary \ref{cor:blockgroups}, every maximal strongly cospectral set in a Cayley graph $X(G,\C)$ is a coset of a subgroup which is then also a maximal strongly cospectral set. We will therefore focus our attention on strongly cospectral sets that contain the identity. If such a set forms a subgroup, we call it a \textsl{strongly cospectral subgroup} and the maximal strongly cospectral set containing the identity will be called the \textsl{maximal strongly cospectral subgroup}.

\section{Normal Cayley graphs}\label{sec:normal}

We now turn to normal Cayley graphs.
We will start by giving some simple necessary conditions for vertices in a normal Cayley graph to be strongly cospectral to the identity. The following is known, but it also follows nicely from the results in the previous chapter.

\begin{lem}\label{lem:centralinvolution}
  Suppose that the vertex $g\in G$ is strongly cospectral to $1$ in the normal Cayley graph $X=X(G,\C)$. Then $g$ has order at most two and lies in the centre of $G$.
\end{lem}
\proof
  We saw before that for each $h\in G$, the map $x\mapsto xh$ is an automorphism of $X$. It is not too hard to see that in a normal Cayley graph, the map $x\mapsto hx$ is also an automorphism.
  Denote the corresponding permutation matrices by $P_h$ and $P'_{h}$, respectively. Then, $P_he_x= e_{xh}$ and $P_h'e_x = e_{hx}$ for all $h,x\in G$.

  Let $E_1,\dots,E_d$ denote the spectral idempotents of the adjacency matrix $A$ of $X$. As before, $P_h$ and $P'_h$ commute with the matrices $E_r$. Since $g$ is strongly cospectral to $1$, we have for all $r$ that $E_re_g=\eps_rE_re_1$, where $\eps_r\in\{\pm1\}$. This implies
  \begin{align*}
    E_re_{g^2} &= P_gE_re_g\\
    &= P_g(\eps_rE_re_1)\\
    & = \eps_rE_re_g\\
    & = \eps_r^2E_re_1\\
    &= E_re_1
  \end{align*}
  for each $r$. Then, since the idempotents sum to the identity, we get
  \begin{align*}
    e_1 &
     = \sum_{r=1}^d E_re_1
     = \sum_{r=1}^d E_re_{g^2}
     = e_{g^2},
  \end{align*}
  and we have shown that $g^2 = 1$.

  Now let $h\in G$ be an arbitrary element. Then $P_hE_re_g = P_h(\eps_rE_re_1)$ implying that
  \[\eps_rE_re_h = E_re_{gh}.\]
  But similarly, we have $ P'_hE_re_g = P'_h(\eps_rE_re_1)$ and so
  \[\eps_rE_re_h = E_re_{hg}.\]
  Thus $E_re_{gh} = E_re_{hg}$ for all $r$, and again by taking the sum over all $r$ we get that $gh=hg$ for all $h$. Therefore, $g$ lies in the centre of $G$.
\qed

Note that the first part of the proof only uses the automorphisms $P_h$ and these are automorphisms even if $X$ is not normal. Thus, it is true for any Cayley graph that a vertex strongly cospectral to the identity is an involution in the group, but it is not necessarily central. Now the following is easy to prove.

\begin{lem}\label{lem:subgroup}
  In a normal Cayley graph $X = X(G,\C)$, the vertices that are strongly cospectral to $1$ form an elementary abelian 2-group and this is a normal subgroup of $G$.
\end{lem}
\proof
  Let $H$ be the set of vertices that are strongly cospectral to $1$ in $X$. By Lemma \ref{lem:blocks}, $H$ is a block under the action of $\Aut(X)$.
  Clearly, every vertex is strongly cospectral to itself, so $1\in H$, and now it follows from Corollary \ref{cor:blockgroups} that $H$ is a subgroup of $G$.
  By Lemma \ref{lem:centralinvolution}, every element of $H$ has order two, so $H$ is an elementary abelian $2$-group. Finally, since the elements of $H$ are central in $G$, it is clear that $H$ is normal in $G$.
\qed

\begin{xmpl}[Non-examples]
  \(\)
  \begin{enumerate}
    \item A Cayley graph of a group of odd order has no strongly cospectral sets, since it has no elements of order two. More generally, in a group of order $2^dm$ where $m$ is odd, a strongly cospectral set has size at most $2^d$.
    \item A normal Cayley graph of the symmetric group on $n$ elements has no strongly cospectral sets, since it has trivial centre.
    \item A normal Cayley graph of a simple group has no strongly cospectral sets by Lemma \ref{lem:subgroup}.
    \item If a Cayley graph of a cyclic group $\Z_n$ has strongly cospectral sets, then $n$ is even and the sets have size two.
    \item Similarly, if there is a strongly cospectral set in a normal Cayley graph of a dihedral group or an extraspecial group, it has size two.
  \end{enumerate}
\end{xmpl}

Lemmas \ref{lem:centralinvolution} and \ref{lem:subgroup} give some good restrictions on the maximal strongly cospectral subgroup in a normal Cayley graph and consequently on the size of a strongly cospectral set in such a graph. We will show later that a strongly cospectral set in a normal Cayley graph can contain at most a third of the vertices (Theorem \ref{thm:n/3}).

\section{Spectrum of a normal Cayley graph}\label{sec:normalspec}

The eigenvalues of normal Cayley graphs and their multiplicities can be calculated using the irreducible characters of the corresponding group. We refer to \cite{james2001reps} for definitions and some basic results from representation and character theory. For details of the following discussion, see \cite[Chapter 11]{godsil2016EKR}.

Let $G$ be a group of order $n$. For an irreducible character $\chi$ of $G$, define a matrix $E_\chi$ with rows and columns indexed by the elements of $G$ by
\[(E_\chi)_{g,h} = \frac{\chi(1)}n\chi(hg^{-1}).\]
It can be verified that $E_\chi$ is idempotent and that if $\psi$ is a character of $G$ different from $\chi$ then $E_\chi E_\psi=0$. Further, if $X=X(G,\C)$ is a normal Cayley graph with adjacency matrix $A$, it can be shown that $AE_\chi=\theta E_\chi$ for some $\theta$. If $\chi_1,\dots,\chi_k$ are all the characters of $G$ such that $AE_{\chi_i}=\theta E_{\chi_i}$, then
\[E_\theta:=\sum_{i=1}^k E_{\chi_i}\]
is the spectral idempotent of $A$ corresponding to the eigenvalue $\theta$. This is the idea behind the proof of the following theorem.
\begin{thm}[{\cite[Theorem 11.12.3]{godsil2016EKR}}]\label{thm:ekr-evals}
  If $X=X(G,\C)$ is a normal Cayley graph and $\chi$ an irreducible character of $G$, then
  \[\theta_\chi =\frac1{\chi(1)}\sum_{c\in \C}\comp{\chi}(c)\]
  is an eigenvalue of $X$ and every eigenvalue can be obtained in this way for some $\chi$.
  Moreover, if $\chi_1,\dots,\chi_k$ are all the irreducible characters such that $\theta_{\chi_i}=\theta$, then $\theta$ has multiplicity
  \[\sum_{i=1}^k\chi_i(1)^2. \qed\]
\end{thm}

We can characterize the vertices that are strongly cospectral to the identity in terms of the irreducible characters of the group. The following was proved by Sin and Sorci in \cite{sin2020}.

\begin{thm}[{\cite[Theorem 2.3]{sin2020}}]\label{thm:characterization}
  Let $X=X(G,\C)$ be a normal Cayley graph. A vertex $g\neq 1$ is strongly cospectral to $1$ if and only if $g$ is a central involution in $G$ and for all irreducible characters $\chi, \psi$ with $\theta_\chi=\theta_\psi$ we have
  \[\frac{\chi(g)}{\chi(1)} = \frac{\psi(g)}{\psi(1)}.  \qed\]
\end{thm}
The proof relies on the fact that an irreducible representation of a group maps a central element to a scalar matrix $cI$. So if $g$ is central and $\chi$ a character corresponding to the irreducible representation $\rho$, then for any $x\in G$ we have
\[\chi(gx) = \tr(\rho(gx)) = \tr(\rho(g)\rho(x)) = \tr(c\rho(x))=c\chi(x),\]
in particular, $\chi(g)=c\chi(1)$. If $g$ is also an involution, then $cI$ has order at most two, so $c=\pm1$ and therefore
\[\frac{\chi(g)}{\chi(1)} =\pm1.\]
We will use this idea again in the proof of Theorem \ref{thm:multbound}.

\section{Multiplicity bound in normal Cayley graphs}\label{sec:multnorm}

We are now ready to prove our first main result which gives an upper bound on the number of pairwise strongly cospectral vertices in a normal Cayley graph in terms of the multiplicities of the eigenvalues of the graph.

\begin{thm}\label{thm:multbound}
  Let \(X=X(G,\C)\) be a normal Cayley graph and let \(H\) be the maximal strongly cospectral subgroup in \(X\). Then, if \(m\) is the multiplicity of some eigenvalue of \(X\) we have
  \[|H|\leq \frac{|G|}m=\frac{|V(X)|}m.\]
\end{thm}
\proof
  Let $\theta$ be an eigenvalue of $X$ and let $\psi_1,\dots,\psi_k$ be a complete set of irreducible characters of $G$ satisfying $\theta_{\psi_i}=\theta$, where $\theta_\psi$ is defined as in Theorem \ref{thm:ekr-evals}. Define $d_i$ to be the degree of $\psi_i$, i.e.\ $d_i:=\psi_i(1)$.
  Then by Theorem \ref{thm:ekr-evals}, the multiplicity of $\theta$ is \(m:=d_1^2+\dots+d_k^2\). Let $\ell$ denote the  index of $H$ in $G$. We will show that $m\leq \ell$.

  We further let $\psi_{k+1},\dots,\psi_n$ be such that $\{\psi_1,\dots,\psi_n\}$ is a complete set of irreducible characters of $G$, and define $d_i$ accordingly for \(i=k+1,\dots,n\). By Theorem \ref{thm:characterization}, we know that for all $i=1,\dots,k$, we have
  \[c_h:=\frac{\psi_i(h)}{d_i} = \frac{\psi_1(h)}{d_1} =\pm 1\quad\text{for all }h\in H.\]

  Let $\rho_i$ denote the irreducible representation corresponding to the character $\psi_i$, for $i=1,\dots,k$. As in the discussion following Theorem \ref{thm:characterization}, we have that $\rho_i(h) =c_hI$ for all $h\in H$.
  Further, since $\rho_i$ is a homomorphism we see that $c_{h_1}c_{h_2} = c_{h_1h_2}$ for all $h_1,h_2\in H$.

  Define a function $\chi:H\to\Cx$ by $h\mapsto c_h$. By the above, this is a homomorphism from $H$ to $\{\pm1\}$, and since $H$ is an abelian group, it is an irreducible character of $H$.
  Further, it is clear by the definition of $\chi$ that for each $i=1,\dots,k$, we have
  \[(\psi_i\downarrow H) = d_i\chi,\]
  where $(\psi\downarrow H)$ denotes the restricted character of $\psi$ to $H$ (defined in the obvious way). If we let $\langle\, ,\,\rangle$ denote the normalized inner product of characters, then this implies that $\langle (\psi_i\downarrow H),\chi\rangle_H = d_i$ for $i=1,\dots,k$.
  We will consider the \textsl{induced character}, $\tilde{\chi} :=\chi\uparrow G$. For convenience, define $\chi':G\to\Cx$ by
  \[\chi'(g) :=
  \begin{cases}
      \chi(g)&\text{if }g\in H\\
      0&\text{otherwise.}
  \end{cases}\]
  Then, the values of the induced character are given by
  \[\tilde{\chi}(g) = \frac1{|H|}\sum_{x\in G}\chi'(x^{-1}gx),\quad\text{for all }g\in G,\]
  in particular, we see that
  \[\tilde{\chi}(1) = \frac{|G|}{|H|} = \ell.\]

  Since $\tilde{\chi}$ is a character of $G$, and $\psi_1, \dots \psi_n$ are the irreducible characters of $G$, we know that $\tilde\chi$ can be written uniquely as
  \[\tilde{\chi} = d'_1\psi_1 +\dots+d'_n\psi_n,\]
  where $d'_i := \langle \tilde{\chi},\psi_i\rangle_G$ are non-negative integers for all $i=1,\dots,n$. By the Frobenius Reciprocity Theorem \cite[Theorem 21.16]{james2001reps}, we have for $i=1,\dots,k$
  \[d'_i = \langle \tilde{\chi},\psi_i\rangle_G = \langle (\psi_i\downarrow H),\chi\rangle_H = d_i,\]
  and so
  \begin{align*}
    \tilde{\chi} &= d_1\psi_1+\dots+d_k\psi_k+d'_{k+1}\psi_{k+1}+\dots+d'_n\psi_n.
  \end{align*}
  Then, evaluating $\tilde\chi$ at $1$, we get
  \begin{align*}
    \ell &=\tilde\chi (1)\\
    & = d_1\psi_1(1)+\dots+d_k\psi_k(1)+d'_{k+1}\psi_{k+1}(1)+\dots+d'_n\psi_n(1)\\
    & = d_1^2+\dots+d_k^2+ d'_{k+1}\psi_{k+1}(1)+\dots+d'_n\psi_n(1)\\
    & = m + K
  \end{align*}
  where \(K\geq 0\) because \(d'_i,\psi_i(1)\geq0\) for all \(i\). It follows that $\ell\geq m$, as required.
\qed

A natural question to ask now is what can we say about the multiplicities of the eigenvalues of a normal Cayley graph? Specifically, what can we say about the largest such multiplicity?
This could be a hard question to answer in general, but we will see in the next sections that for certain groups, the Cayley graphs tend to have some eigenvalues with large multiplicities.

\section{Cubelike graphs}
Recall that a cubelike graph is a Cayley graph of an elementary abelian 2-group. Such a graph is normal, because the group is abelian. Throughout the rest of this paper, we will be considering abelian groups. We will think of them as additive and the identity will be called zero.

Cubelike graphs are in many ways a natural choice when looking for large strongly cospectral sets. Firstly, most cubelike graphs (in some loose sense of the word ``most'') have at least a pair of strongly cospectral vertices. This follows from a result of Bernasconi et al \cite{bernasconi}, and the fact mentioned earlier that perfect state transfer implies strong cospectrality.
\begin{thm}\label{thm:cubelikepst}{\cite[Theorem 1]{bernasconi}}
  A cubelike graph $X(\Z_2^d,\C)$ satisfying
  \[\sigma:=\sum_{c\in\C}c\neq 0\]
  has perfect state transfer from $0$ to $\sigma$ at time $\pi/2$. \qed
\end{thm}
Secondly, all the elements of $\Z_2^d$ are central involutions, and so every vertex in a cubelike graph is a candidate for being strongly cospectral to $0$. Therefore, one might think that a cubelike graph could have large strongly cospectral sets. We will however show that the maximal strongly cospectral subgroup in a Cayley graph of $\Z_2^d$ has dimension less than $d/2$.

For abelian groups, every irreducible character has degree one and is therefore equal to the corresponding representation.
It is left as a fun exercise for the reader to verify that the characters of an abelian group, $G$, are eigenvectors of the Cayley graph $X(G,\C)$ (for any connection set $\C$) and the eigenvalue for the eigenvector $\psi$ is
\[\psi(\C) := \sum_{c\in \C} \psi(c).\]

Since the elements of an elementary abelian 2-group have order two, every character takes values in $\{\pm1\}$, and so we see that every eigenvalue of a cubelike graph is an integer with the same parity as $n:=|\C|$ and lies on the interval $\left[-n,n\right]$. It turns out that the eigenvalues of a cubelike graph are, to put it very vaguely, close to being normally distributed about the origin. That is, the eigenvalues that are close to $n$ and $-n$ have small multiplicities whereas there will be some eigenvalues close to $0$ with large multiplicities.

Let $X=X(\Z_2^d,\C)$ be a cubelike graph with degree $n=|\C|$ and let $A$ be its adjacency matrix. Each eigenvalue of $X$ can then be written on the form $n-2r$, with $r=0,\dots,n$. Let $m_r$ be the multiplicity of the eigenvalue $n-2r$, with $m_r=0$ if $n-2r$ is not an eigenvalue. Then, $m_0+\dots+m_n=|V(X)|=2^d$.
Consider the matrix $A^2$. Since $X$ is regular with degree $n$, every diagonal entry of $A^2$ is $n$, and so it has trace $2^dn$. Further, the eigenvalues of this matrix are the squares of the eigenvalues of $A$ and now, since the sum of the eigenvalues of a matrix is equal to its trace, we have shown the following.

\begin{lem}\label{lem:cubeid}
  Let $X$ be a cubelike graph on $2^d$ vertices, with degree $n$ and for $r=0,\dots,n$, let $m_r$ be the multiplicity of the eigenvalue $n-2r$ of $X$. Then
  \[\sum_{r=0}^n(n-2r)^2m_r = 2^dn.  \qed\]
\end{lem}
We shall use this identity to prove Theorem \ref{thm:cubemult}.

\section{A bound on strongly cospectral sets in cubelike graphs} \label{sec:cubelikebound}

We will give a lower bound on the largest multiplicity of an eigenvalue in a cubelike graph. Combining it with Theorem \ref{thm:multbound} directly gives an upper bound on the size of a strongly cospectral set in a cubelike graph in terms of the dimension of the group.

Before we do that we need one more lemma about the connection between the spectrum of a graph and its complement. We will omit the proof here, but the idea is that a regular graph and its complement have the same eigenvectors (see for example \cite[Lemma 8.5.1]{godsilAGT}), and that the degree is an eigenvalue with multiplicity equal to the number of components.

\begin{lem}\label{lem:complements}
  Let $X$ be a connected, vertex-transitive graph and denote by $\comp X$ its complement. If $\theta$ is an eigenvalue of $\comp X$ different from the degree of $\comp X$, with multiplicity $m$, then $X$ has an eigenvalue with multiplicity $m$.\qed
\end{lem}

We can now prove the main result of this section.

\begin{thm}\label{thm:cubemult}
  Let $X=X(\Z_2^d,\C)$ be a cubelike graph, with $d\geq3$ and let $q=\frac d2$. Then $X$ has an eigenvalue with multiplicity larger than $2^q$.
\end{thm}
  Note that the theorem does not hold for \(d=2\); a \(4\)-cycle does not have an eigenvalue with multiplicity larger than two.

\proof
  Suppose first that \(X\) is not connected. Then each connected component of \(X\) is a cubelike graph on fewer than \(2^d\) vertices and they are all isomorphic. Suppose \(X\) has \(c\) components, let \(Y\) be one of them and let \(2^{d'}\) be its number of vertices. Then \(2^{d'}c=2^d\), so \(c=2^{d-d'}\). Further, each eigenvalue of \(Y\) with multiplicity \(m\) is an eigenvalue of \(X\) with multiplicity \(cm\).

  If \(d'\geq3\), we may assume inductively, that \(Y\) has an eigenvalue with multiplicity \(m>2^{d'/2}\), and so \(X\) has an eigenvalue with multiplicity
  \[cm>2^{d-d'}\cdot2^{d'/2} = 2^{d-d'/2}>2^{d-d/2}=2^q.\]
  If \(d'\in \{0,1\}\), an eigenvalue of \(Y\) with multiplicity one gives an eigenvalue of \(X\) with multiplicity \(c=2^{d-d'}>2^q\). Finally if \(d'=2\), then \(Y\in\{C_4,K_4\}\) thus it has an eigenvalue with multiplicity at least two, and so \(X\) has an eigenvalue with multiplicity at least
  \[2c= 2\cdot 2^{d-2}=2^{d-1}>2^q.\]

  Then suppose \(X\) is connected. As before, let \(n:=|\C|\) and let \(m_r\) denote the multiplicity of the eigenvalue \(n-2r\). Then
  \[\sum_{r=0}^nm_r =2^d\]
  and since $X$ is a connected, $n$-regular graph, $n$ is an eigenvalue with multiplicity one and $-n$ has multiplicity at most one. Therefore, \(m_0=1\) and \(m_n\in\{0,1\}\).
  We will assume by way of contradiction that \(m_r\leq 2^q\) for all \(r\). Then the above gives
  \[2^{2q}=\sum_{r=0}^n m_r \leq 2+2^q(n-1),\]
  which implies
  $n\geq 2^q- 2/2^q+1>2^q,$
  (since \(q> 1\)).

  Assume first that both $d$ and \(n\) are even, so \(q\) is an integer and let \(k:=\frac n2\). Recall the identity from Lemma \ref{lem:cubeid},
  \begin{align*}
    2^dn = \sum_{r=0}^n(n-2r)^2m_r.
  \end{align*}
  We will see that under the assumption that $m_r\leq 2^q$ for all $r$, we can use this to derive a quadratic inequality, $n^2+an+b\leq 0$, with negative discriminant, yielding a contradiction since $n$ is an integer.
  Notice that
  \[\sum_{r=0}^n(n-2r)^2m_r = n^2(1+m_n) +\sum_{r=1}^{k-1}(n-2r)^2(m_r+m_{n-r}).\]
  We will split this sum into two parts. Define $t:=k-2^{q-1}$. We see that $0<t<k-1,$ and now Lemma \ref{lem:cubeid} gives
  \begin{align}
    2^dn-n^2 &= n^2m_n + \sum_{r=1}^{k-1} (n-2r)^2(m_r+m_{n-r})\nonumber\\[0.1cm]
    & = n^2m_n + \sum_{r=1}^{t}(n-2r)^2(m_r+m_{n-r})\nonumber\\[0.1cm]
    &\phantom{==} + \sum_{r=t+1}^{k-1}(n-2r)^2(m_r+m_{n-r})\nonumber\\[0.1cm]
    & \geq (n-2t)^2 \left(m_n+\sum_{r=1}^{t}(m_r+m_{n-r})\right)\label{eq:RHS1}\\[0.1cm]
    &\phantom{==} + \sum_{r=t+1}^{k-1}(n-2r)^2(m_r+m_{n-r}).\label{eq:RHS2}
  \end{align}
  We can rewrite (\ref{eq:RHS2}) as
  \begin{align*}
    \sum_{r=t+1}^{k-1}(n-2r)^2(m_r+m_{n-r}) &= \sum_{r=1}^{2^{q-1}-1}(2r)^2(m_{k-r}+m_{k+r}).
  \end{align*}
  and using the fact that $m_1+\dots+m_n = 2^d-1$, we see that (\ref{eq:RHS1}) is equal to
  \begin{align*}
    &\phantom{==}(n-(n-2^{q}))^2\left(2^d-1-m_k-\sum_{r = t+1}^{k-1} (m_r+m_{n-r})\right)\\[0.1cm]
    & = 2^d(2^d-1-m_k)-\sum_{r =1}^{2^{q-1}-1}2^d(m_{k-r}+m_{k+r}).
  \end{align*}
  Putting these together, and assuming that \(m_r\leq 2^q\) for all \(r\), we get
  \begin{align*}
    &\phantom{==}2^dn-n^2\\
    &\geq 2^d(2^d-1-m_k) -\sum_{r =1}^{2^{q-1}-1}\left(2^d-(2r)^2\right)(m_{k-r}+m_{k+r})\\
    & \geq 2^{2d}-2^d-2^{d+q}-2^{q+1}\sum_{r=1}^{2^{q-1}-1}\left(2^{2q}-(2r)^2\right)\\[0.1cm]
    & = 2^{2d}-2^d-2^{d+q}-2^{q+1}\cdot\frac{2^{q-1}\left(2^{2q+1}-3\cdot2^q-2\right)}3\\[0.1cm]
    & = 2^{2d}-2^d+\frac{2^{d+1}-2^{2d+1}}3
  \end{align*}
  and so
  \begin{align}
  n^2-2^dn+2^{2d}-2^d+ \frac{2^{d+1}-2^{2d+1}}3\leq 0.
  \label{eq:wooo}
\end{align}
  This is a quadratic inequality in \(n\) with discriminant
  \begin{align*}
    &\phantom{==}2^{2d}-4\left(2^{2d}-2^d+\frac{2^{d+1}-2^{2d+1}}3\right)\\[0.1cm]
    & = -3\cdot 2^{2d}+4\cdot 2^d-\frac{8\left(2^{d}-2^{2d}\right)}3\\[0.1cm]
    & = \frac13 \left(2^{2d}(8-9)+2^{d+2}(3-2)\right)\\[0.1cm]
    & = \frac13(2^{d+2}-2^{2d})\\[0.1cm]
    & < 0,
  \end{align*}
  since \(d>2\). Therefore, (\ref{eq:wooo}) never holds when \(d\geq 3\), and we have reached the required contradiction.

  Next we consider the case where $n$ is even and $d$ is odd. Here, $q$ is no longer an integer so we need to split our sum differently. We let $q_0:=(d-1)/2$ and define $t=k-2^{q_0-1}$. We use the same technique as before to get
  \begin{align*}
    n^2-2^dn +2^{2d-1}-2^{d-1}+\frac{2^{q+q_0+1}-2^{d+q+q_0}}3\leq 0,
  \end{align*}
  which has discriminant
  \begin{align*}
    &\phantom{==}2^{2d}-4\left(2^{2d-1}-2^{d-1}+\frac{2^{q+q_0+1}-2^{d+q+q_0}}3\right)\\
    & = \frac13\left(-3\cdot2^{2d}+3\cdot2^{d+1}+4\cdot2^{2d-\frac12}-4\cdot2^{d+\frac12}\right)\\[0.1cm]
    & = \frac13\left(2^{2d}(2\sqrt{2}-3)+2^{d+1}(3-2\sqrt2)\right)\\[0.1cm]
    & =\frac{3-2\sqrt2}3\left(2^{d+1}-2^{2d}\right)\\[0.1cm]
    & < 0,
  \end{align*}
  since \(d>2\), so again we have a contradiction. Note that if $d=3$ then $t=k-1$, so the sum in (\ref{eq:RHS2}) is empty, but the same argument still holds.

  Now assume that \(n\) is odd. Then the complement, $\comp X$, of $X$ is a cubelike graph on $2^d$ vertices with even degree, so by the above it has an eigenvalue with multiplicity $m > 2^q$. If this eigenvalue is different from the degree of $\comp X$ then by Lemma \ref{lem:complements}, $X$ also has an eigenvalue with multiplicity $m$.

  So suppose that the degree $k$ of $\comp X$ is an eigenvalue with multiplicity $m$.
  The components of $\comp X$ are isomorphic, connected cubelike graphs of degree $k$. If $Y$ is one such component and $\theta$ is an eigenvalue of $Y$ with multiplicity $m'$ then $\theta$ is an eigenvalue of $\comp X$ with multiplicity $cm'$ where $c$ is the number of components. But $k$ is an eigenvalue of $Y$ with multiplicity one, and so we must have $c=m$. Further, if $Y$ is not the one-vertex graph, then it has another eigenvalue, $\theta\neq k$ and this is an eigenvalue of $\comp X$ with multiplicity $m'\geq m$.
  As before, this gives an eigenvalue of $X$ with multiplicity $m'>2^q$ by Lemma \ref{lem:complements}.

  We are left with the case where $\comp X$ is edgeless, but then $X$ is complete and has eigenvalue $-1$ with multiplicity $2^d-1>2^q$ and this completes the proof.
\qed

Now, combining Theorem \ref{thm:cubemult} and Theorem \ref{thm:multbound} we get the following upper bound on the size of a strongly cospectral set in a cubelike graph.

\begin{cor}\label{cor:cubebound}
  In a cubelike graph on \(2^d\) vertices, with \(d\geq 3\), a strongly cospectral set has size at most
  \(2^{\lceil d/2\rceil-1}\).
\end{cor}
\proof
  Let \(X\) be a cubelike graph on \(2^d\) vertices with \(d\geq3\), let $s$ be size of a maximal strongly cospectral set in $X$ and let \(m\) be the largest multiplicity of an eigenvalue of \(X\). Then by the above theorem, $m>2^{d/2}$ and so Theorem \ref{thm:multbound} gives
  \[s\leq \frac{2^d}{m}< 2^{d/2}.\]
  But by Lemma \ref{lem:subgroup}, $s$ is a power of two, so we get \(s\leq2^{\lceil d/2\rceil-1}\).
\qed

\begin{xmpl}\label{ex:tight}$ $
  \begin{enumerate}
    \item For $d\in\{3,4\}$, we get $2^{\lceil d/2\rceil-1}=2$. In every hypercube, antipodal vertices are strongly cospectral, so the $3$-cube and the $4$-cube are examples of the bound in Corollary \ref{cor:cubebound} being tight.
    \item Let $d=5$, then $2^{\lceil d/2\rceil-1}=4$. There are exactly twelve cubelike graphs (in six complementary pairs) on $2^5=32$ vertices that have strongly cospectral sets of size four. Three of them are shown in Figure \ref{fig:ex32}. The connection sets of all six graphs are given in Appendix A.
    \begin{figure}[h!]
      \centering
      \includegraphics[scale=0.35]{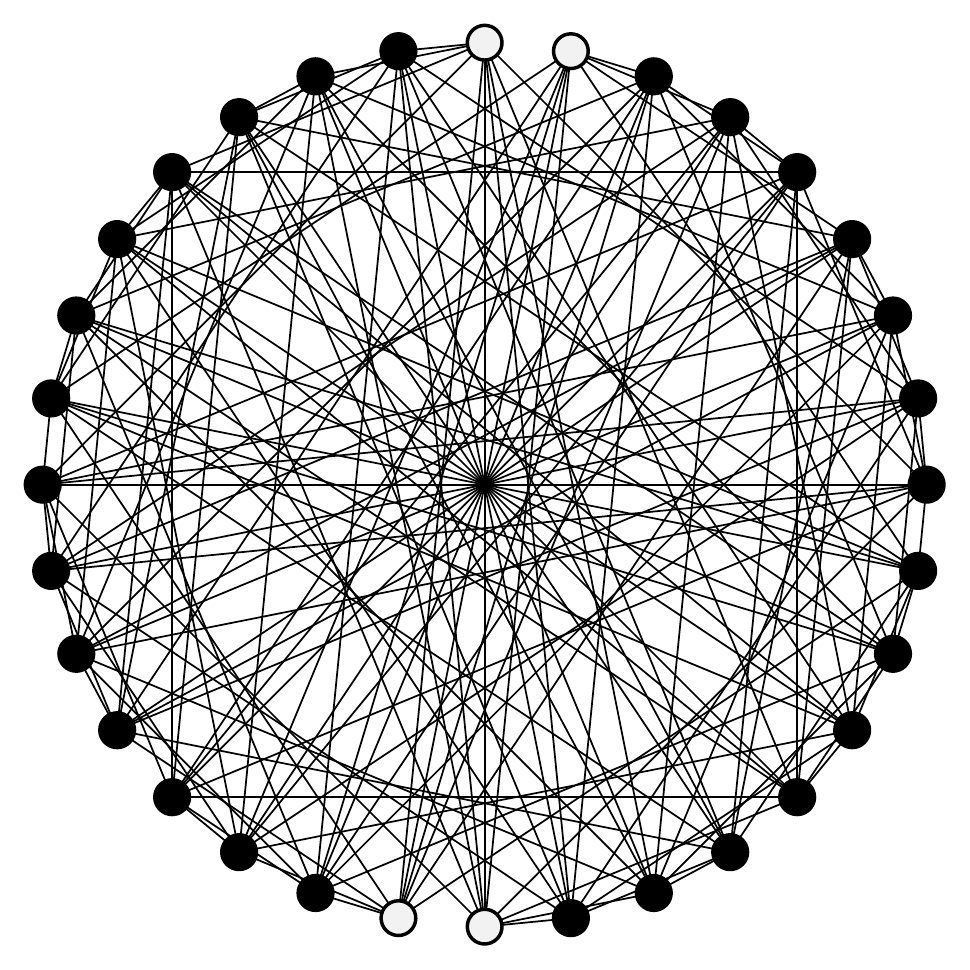}\hspace{0.2cm}
      \includegraphics[scale=0.35]{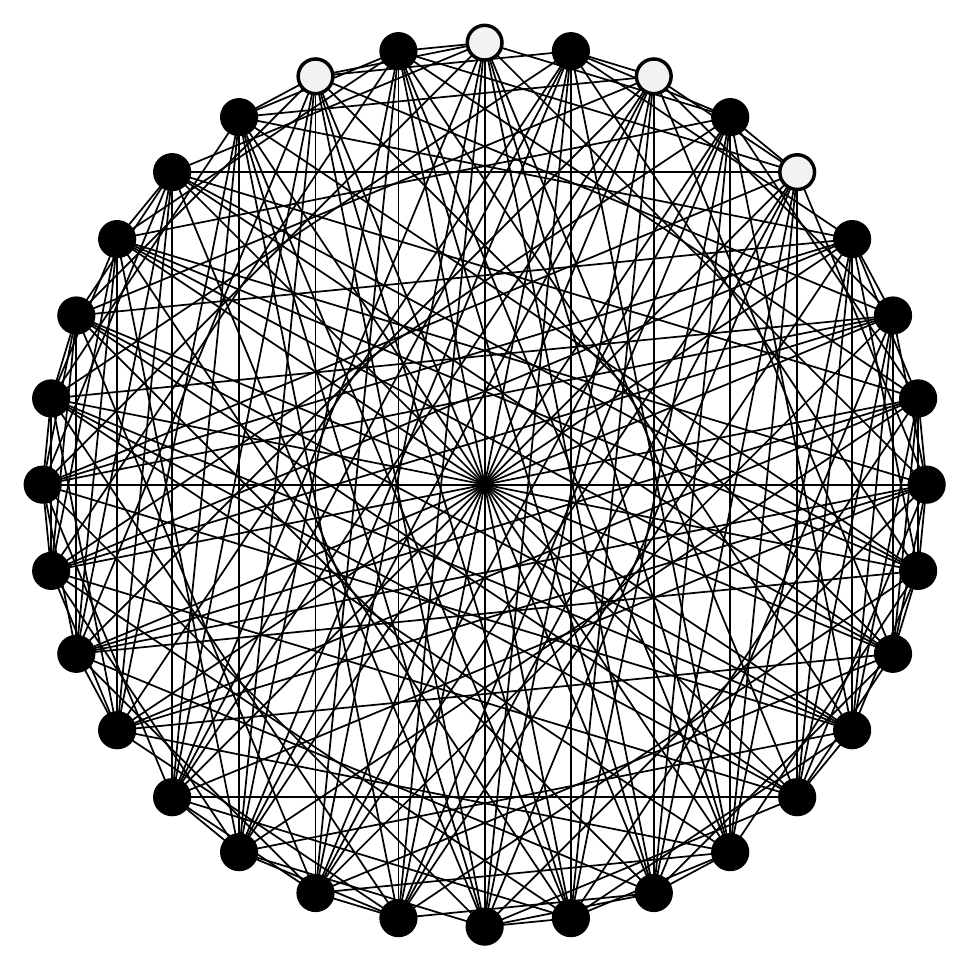}\hspace{0.2cm}
      \includegraphics[scale=0.35]{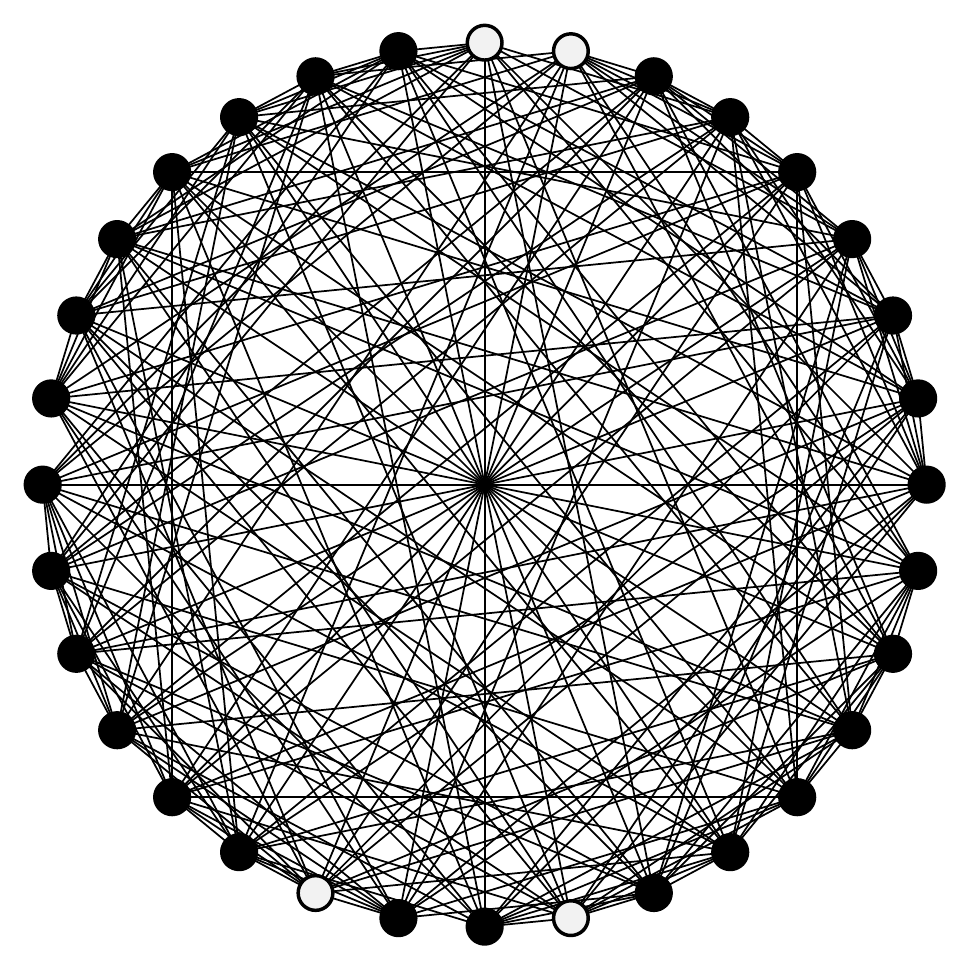}
      \caption{Cubelike graphs on $32$ vertices of degree $10, 13$ and $14$, respectively. The vertices in $H$ are white.}
      \label{fig:ex32}
    \end{figure}
    \item For $d=6$ we can again have at most a strongly cospectral set of size four. There are many such examples, one of which is shown in Figure \ref{fig:ex64}
    \begin{figure}[h!]
      \centering
      \includegraphics[scale=0.75]{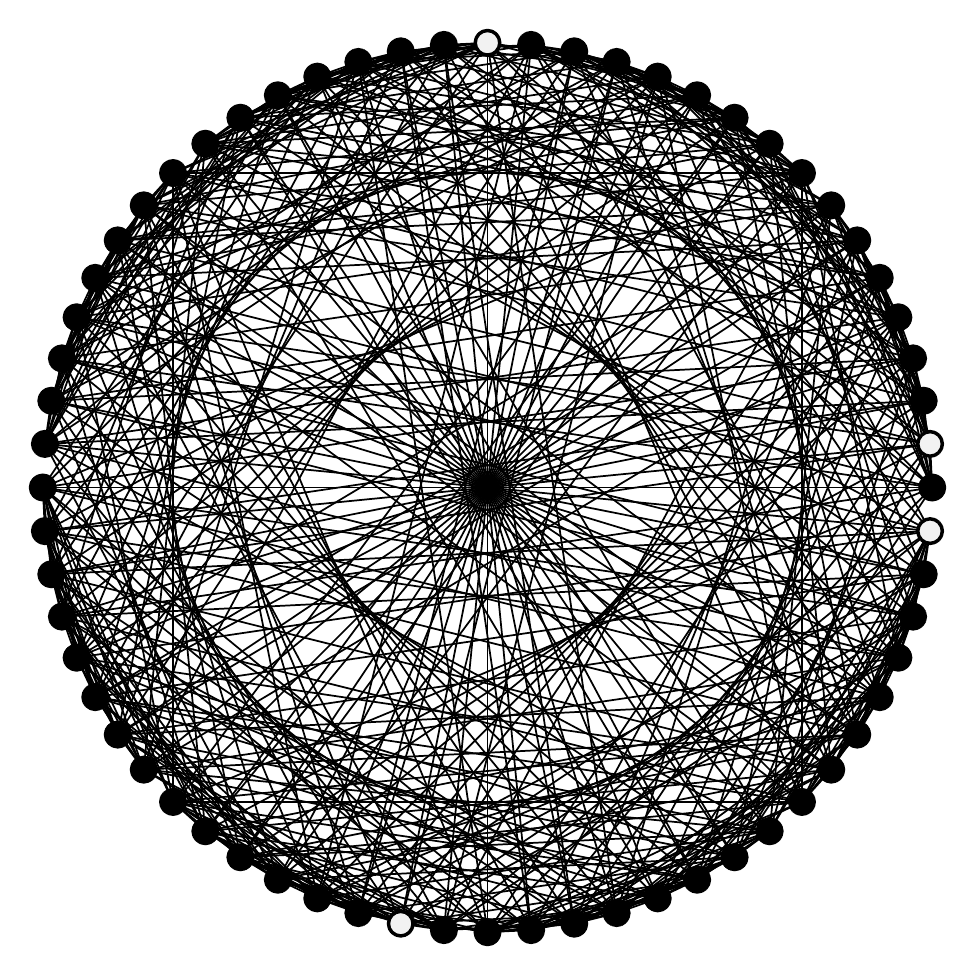}
      \caption{A cubelike graph on $64$ vertices of degree $16$. The vertices in $H$ are white.}
      \label{fig:ex64}
    \end{figure}
  \end{enumerate}
\end{xmpl}
Example \ref{ex:tight} shows that the bound given in Corollary \ref{cor:cubebound} is tight for dimensions at most six. We do not know whether the bound is tight in general and in fact we currently have no examples of cubelike graphs (or even vertex-transitive graphs) with more than four pairwise strongly cospectral vertices.

\section{Cayley graphs of \(\Z_2^{d_1}\times\Z_4^{d_2}\)}

Notice that the proof of Theorem \ref{thm:cubemult}, relies on $X$ being cubelike, only in that its eigenvalues may be written on the form $n-2r$, where $n$ is the degree and $r=0,\dots,n$. We claim that this is true for Cayley graphs of all groups of the form $\Z_2^{d_1}\times\Z_4^{d_2}$.

Let $X=X(\Z_4^d,\C)$ be a Cayley graph and let $\psi$ be a character of $\Z_4^d$. Then $\psi$ is a homomorphism of $\Z_4^d$ and takes values in $\{\pm 1,\pm i\}$. Further, $\psi$ is an eigenvector with eigenvalue
\[\psi(\C) = \sum_{c\in\C}\psi(c).\]
Recall that by definition, $\C$ is inverse closed. This means that whenever $\psi(c)=i$ appears in the sum, we have $\psi(c^{-1})=-i$ also appearing in the sum, so they cancel out. So, $\psi(\C)$ is again a sum of ones and negative ones, and the number of terms again has the same parity as $n:=|\C|$. Therefore, every eigenvalue can be written on the form $n-2r$ with $r= 0,\dots, n$. Clearly, since this is true for Cayley graphs for the groups $\Z_2^d$ and $\Z_4^d$, it holds in general for groups of the form $\Z_2^{d_1}\times\Z_4^{d_2}$.

We can therefore recycle the proof of Theorem \ref{thm:cubemult}, taking $d:=d_1+2d_2$, to prove the following.

\begin{thm}\label{thm:Z4dmult}
  Let $X=X(\Z_2^{d_1}\times\Z_4^{d_2},\C)$ be a Cayley graph, with $d:=d_1+2d_2\geq 3$ and define $q:=\frac d2$. Then $X$ has an eigenvalue with multiplicity larger than $2^q$.\qed
\end{thm}

It is then easy to prove a more general version of Corollary \ref{cor:cubebound}, but a proof can also be found in \cite[Theorem 5.7.2]{ArnadottirThesis}.

\begin{thm}\label{thm:Z2Z4}
  Let $d_1, d_2 \geq 0$ be integers such that $d:= d_1 + 2d_2 \geq 3$. Then a strongly cospectral subgroup of $\Z_2^{d_1}\times\Z_4^{d_2}$ has order at most $2^{\lceil d/2\rceil -1}$.\qed
\end{thm}

\begin{cor}\label{cor:Z4dbound}
  Any strongly cospectral set in a Cayley graph $X(\Z_4^d,\C)$, where \(d\geq 2\), has size at most
  \(2^{d-1}\).\qed
\end{cor}

We can now prove the following general bound on the size of a strongly cospectral set in a normal Cayley graph.

\begin{thm}\label{thm:n/3}
  In a normal Cayley graph, $X=X(G,\C)$ on at least five vertices, a strongly cospectral set has size at most $\frac{|V(X)|}{3}$.
\end{thm}
\proof
  Let $H$ be the maximal strongly cospectral subgroup. If $H=G$, all the vertices are pairwise strongly cospectral, but this is impossible by Godsil and Smith \cite[Lemma 10.1]{godsil2017strongly}.

  Then supposing by way of contradiction that $|H|>|G| / 3$, the only possibility is that $|G:H|=2$. Since $H$ is contained in the centre, $Z(G)$, this means that the quotient group $G/Z(G)$ has order at most two and is therefore cyclic. It is a known fact from group theory that in this case, $G$ has to be abelian. Then $G$ is an abelian group containing an elementary abelian 2-subgroup $H$ with index two, and clearly it follows that $G=\Z_2^{d_1}\times\Z_4^{d_2}$, where $d_2\in\{0,1\}$ and we apply Theorem \ref{thm:Z2Z4}.\qed

\begin{note*}
  It can be shown that a Cayley graph for $\Z_{2}^{d_1}\times \Z_4^{d_2}$ is isomorphic to a cubelike graph for $\Z_2^d$ with $d=d_1+2d_2$. The results of this section can also be concluded from this fact.
\end{note*}

\section{Sets of size four}

In this section, we show through construction that for all $d\geq 5$, there exists a connected cubelike graph on $2^d$ vertices that has strongly cospectral sets of size four. %We do this by constructing such graphs.
We then use this to show that for any positive integer $n$ divisible by $32$, there exists a Cayley graph of an abelian group on $n$ vertices with strongly cospectral sets of size four.
First, we refine Theorem \ref{thm:ekr-evals} for abelian groups:
\begin{thm}\label{thm:abel-evals}
  Let \(G\) be an abelian group. The vertices \(0\) and \(c\) are strongly cospectral in \(X(G,\C)\) if and only if \(2c = 0\) and for any two characters \(\phi\) and \(\psi\), if \(\phi(\C) = \psi(\C)\) then \(\phi(c) = \psi(c)\).\qed
\end{thm}

\begin{figure}[h!]
  \centering
  \includegraphics[scale=0.55]{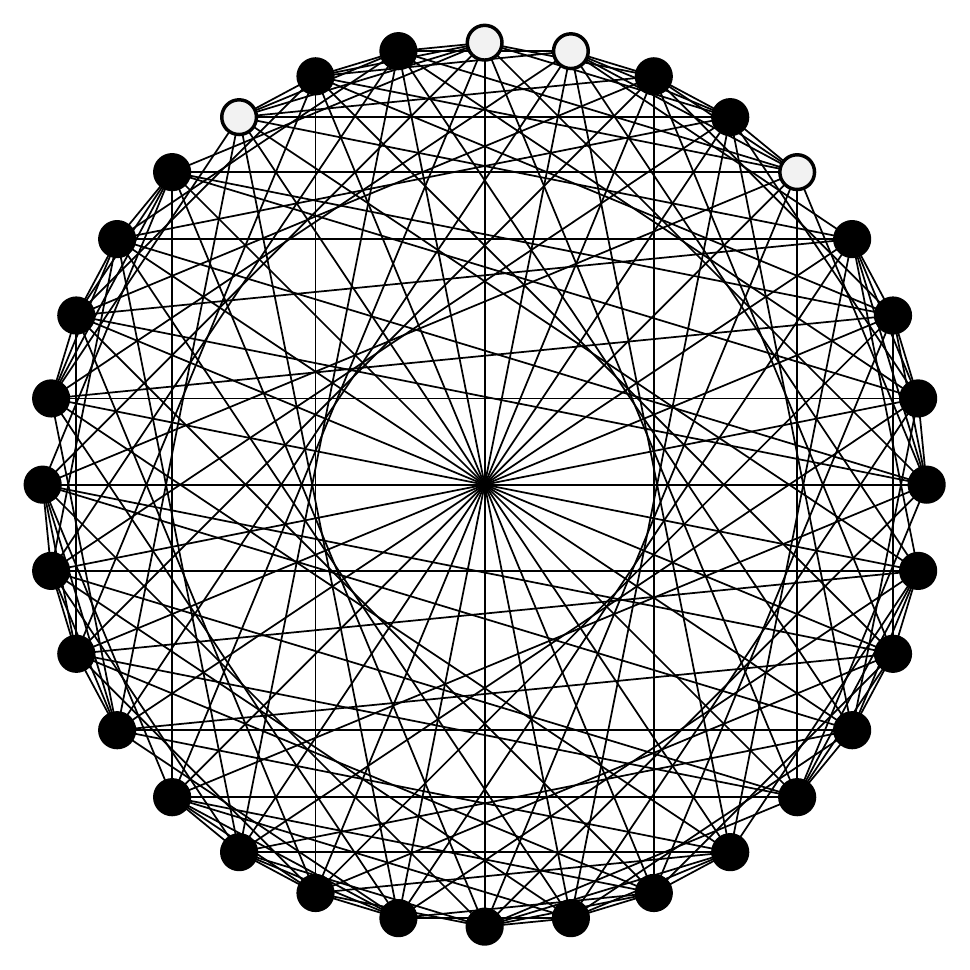}
  \caption{The graph $X$ defined above with $d=5$.}
  \label{fig:alld5}
\end{figure}

The constructions are different for odd and even dimensions.
Let $d$ be an odd integer with $d\geq 5$. Denote by $e_1,\dots,e_d$ the standard basis vectors of $G:=\Z_2^d$.
Define the sets
\begin{align*}
  C_1&:=\{e_1,\,e_2,\dots,\,e_d\},\\
  C_2&:=\{e_1+e_2,\,e_1+e_3,\,e_2+e_3\},\\
  C_3&:=\{e_i+e_j: 4\leq i<j\leq d\},\\
  C_4&:=\{e_1+e_2+e_3+e_i:4\leq i\leq d\},
\end{align*}
let
\(\C:=C_1\cup C_2\cup C_3\cup C_4\), and let $X:=X(G,\C)$. Then
\[\sigma:=\sum_{c\in \C}c = e_1+\cdots +e_d\] and
$X$ has degree
\(2d+(d^2-7d+12)/2.\)
Figure \ref{fig:alld5} shows this graph in dimension $5$. %It has degree $11$.

\begin{thm}\label{thm:exodd}
  Let $X$ be the cubelike graph described above for some odd dimension $d\geq5$ and let $H\leq \Z_2^d$ be the maximal strongly cospectral subgroup in $X$. Then $|H|\geq 4$.
\end{thm}
\proof
  By Theorem \ref{thm:cubelikepst},
  there is perfect state transfer between 0 and \(\sigma = e_1+\dots+e_d\), thus they are strongly cospectral.
  We will use Lemma \ref{thm:abel-evals} to show that the vertex \(g :=e_1+e_2+e_3\) is strongly cospectral to \(0\). Note that since \(d\geq 5\), we have $g\neq \sigma$ and so this implies that $\langle g,\sigma\rangle$ is a subgroup of order four, contained in the subgroup $H$, of elements that are strongly cospectral to zero.

  Since every element of \(G\) has order two, it suffices to show that if \(\psi,\phi\) are characters of \(G\) with $\psi(\C) = \phi(\C)$
  then \(\psi(e_1+e_2+e_3) = \phi(e_1+e_2+e_3)\).
  Consider an arbitrary character \(\chi\) of \(G\). Since $G$ is an abelian group, $\chi:G\to\Cx^*$ is a homomorphism. For convenience, define $\eta_\chi:=\chi(e_1)\chi(e_2)+\chi(e_1)\chi(e_3) + \chi(e_2)\chi(e_3)$. We have
  \begin{align}\label{eq:chiC1}
    \chi(\C)&=\chi(C_1)+\chi(C_2)+\chi(C_3)+\chi(C_4)\\
    \nonumber &=\sum_{i=1}^d\chi(e_i)+\chi(e_1+e_2)+\chi(e_1+e_3) + \chi(e_2+e_3)\\
    \nonumber & + \sum_{4\leq i< j\leq d}\hspace{-0.3cm}\chi(e_i+e_j)+\sum_{i=4}^d \chi(e_1+e_2+e_3+e_i)\\
    \nonumber &= \sum_{i=1}^d\chi(e_i)+\eta_\chi
    + \sum_{4\leq i< j\leq d}\hspace{-0.3cm} \chi(e_i)\chi(e_j) + \chi(e_1)\chi(e_2)\chi(e_3)\sum_{i=4}^d\chi(e_i).
  \end{align}
  Since every element of \(G\) has order two, we know that \(\chi(g)=\pm1\) for all \(g\in G\). For further convenience, let \(d':=d-3\) and consider the set
  \(C_1':=\{e_4,\dots,e_d\}.\)
  Let \(p_\chi\) denote the number of elements \(c\in C_1'\) such that \(\chi(c)=1\) and \(n_\chi\) the number of elements \(c\in C_1'\) such that \(\chi(c)=-1\), i.e.
  \[p_\chi = |\chi^{-1}(1)\cap C_1'|\quad\text{and}\quad n_\chi = |\chi^{-1}(-1)\cap C_1'|.\]
  Then \(p_\chi +n_\chi = d'\) and we have
  \[\sum_{i=4}^d\chi(e_i) = p_\chi-n_\chi = 2p_\chi-d'.\]
  Furthermore,
  \begin{align*}
    \sum_{4\leq i< j\leq d}\hspace{-0.3cm}\chi(e_i)\chi(e_j) & = \binom{p_\chi}2+\binom{n_\chi}2-p_\chi n_\chi\\
    & =\frac12\big(p_\chi(p_\chi-1)+(d'-p_\chi)(d'-1-p_\chi)-2p_\chi(d'-p_\chi)\big)\\
    & =\frac12\big(4p_\chi^2-4d'p_\chi+d'(d'-1)\big)\\
    & = 2p_\chi^2-2d'p_\chi+\frac{d'(d'-1)}2.
  \end{align*}
  Putting this together with Equation \ref{eq:chiC1}, we obtain
  \begin{align}\label{eq:chiC2}
    \chi(\C) &=\chi(e_1)+\chi(e_2)+\chi(e_3) +2p_\chi-d'+\eta_\chi\\[0.1cm]
    \nonumber & + 2p^2_\chi -2d'p_\chi +\frac{d'(d'-1)}2 \pm(2p_\chi-d')
  \end{align}
  where the $\pm$ depends on the value of $\chi(e_1+e_2+e_3)$.
  Now let \(\psi, \phi\) be characters of \(G\) such that \(\psi(\C)=\phi(\C)\) and suppose by way of contradiction that \(\psi(e_1+e_2+e_3)\neq \phi(e_1+e_2+e_3)\). We may assume without loss of generality that \(\psi(e_1+e_2+e_3)=1\) and \(\phi(e_1+e_2+e_3)=-1\). Define \(p_\psi\) and \(p_\phi\) as before.
  Since
  \[\phi(e_1)\phi(e_2)\phi(e_3)=\phi(e_1+e_2+e_3) = -1,\]
  we have two possibilities: either \(\phi(e_1)=\phi(e_2)=\phi(e_3) = -1\) or exactly one out of the two is \(-1\) and the other two are \(1\). We can see that in both cases, we have
  \[\phi(e_1)+\phi(e_2)+\phi(e_3) +\phi(e_1)\phi(e_2) + \phi(e_1)\phi(e_3)+\phi(e_2)\phi(e_3) = 0.\]
  Similarly, we have two cases for \(\psi(e_1),\psi(e_2),\psi(e_3)\): either they are all one, or exactly two of them are \(-1\). It follows that
  \[\psi(e_1)+\psi(e_2)+\psi(e_3) +\psi(e_1)\psi(e_2) + \psi(e_1)\psi(e_3) + \psi(e_2)\psi(e_3) \in\{-2,6\}.\]
  Combining this with Equation \ref{eq:chiC2} we get
  \[\phi(\C) = 2p^2_\phi -2d'p_\phi +\frac{d'(d'-1)}2\]
  and
  \[\psi(\C) =
  \begin{cases}
    2p^2_\psi -2d'p_\psi +\frac{d'(d'-1)}2 +2(2p_\psi-d')-2,& \text{or}\\[0.2cm]
    2p^2_\psi -2d'p_\psi +\frac{d'(d'-1)}2 +2(2p_\psi-d') +6
  \end{cases}
  \]
  Now \(\psi(\C)=\phi(\C)\) implies
  \[2p_\phi^2-2d'p_\phi - (2p_\psi^2-2d'p_\psi+2(2p_\psi-d'))\in\{-2,6\}\]
  and so
  \[p_\phi^2-d'p_\phi - 2p_\psi+d'-p_\psi^2+d'p_\psi\in\{-1,3\}.\]
  Recall that \(d\) is odd, so \(d'=d-3\) is even; let \(d'=2z\) with \(z\in \Z\). Further, since \(\sigma\) is strongly cospectral to 0, we know that \(\phi(\sigma)=\psi(\sigma)\) and it follows that \(p_\psi\) and \(p_\phi\) have different parity. Suppose first \(p_\psi=2x\) and \(p_\phi=2y+1\), with \(x,y\in\Z\). Then
  \begin{align*}
    &(2y+1)^2-2z(2y+1) - 4x+2z-4x^2+4zx\\
    & = 4y^2+4y+1-4zy-4x+4x^2+4zx \in\{-1,3\}
  \end{align*}
  but this implies
  \[4(y^2+y-zy-x+x^2+zx) \in\{-2,2\}\]
  which is impossible since $x,y,z$ are integers. Then suppose \(p_\psi =2x+1\) and \(p_\phi=2y\). Then
  \begin{align*}
    & 4y^2-4zy-2(2x+1)+2z-(2x+1)^2+2z(2x+1)\\
    & = 4y^2-4zy-4x-2+2z-4x^2-4x-1+4zx+2z\\
    & = 4y^2-4zy-8x +4z-3\in\{-1,3\}
  \end{align*}
  which implies
  \[4(y^2-zy-2x +z)\in\{2,6\},\]
  again impossible.

  We conclude that whenever \(d\) is odd and \(\phi\) and \(\psi\) are characters of \(G\) such that \(\phi(\C)=\psi(\C)\), then \(\phi(e_1+e_2+e_3)=\psi(e_1+e_2+e_3)\) and so \(e_1+e_2+e_3\) is strongly cospectral to zero in $X$. Therefore,
  \[\{0,\,g,\,\sigma,\,g+\sigma\}\subseteq H\]
  and so $|H|\geq 4$.
\qed

\begin{figure}[h!]
  \centering
  %\vspace{-0.3cm}
  \includegraphics[scale=0.75]{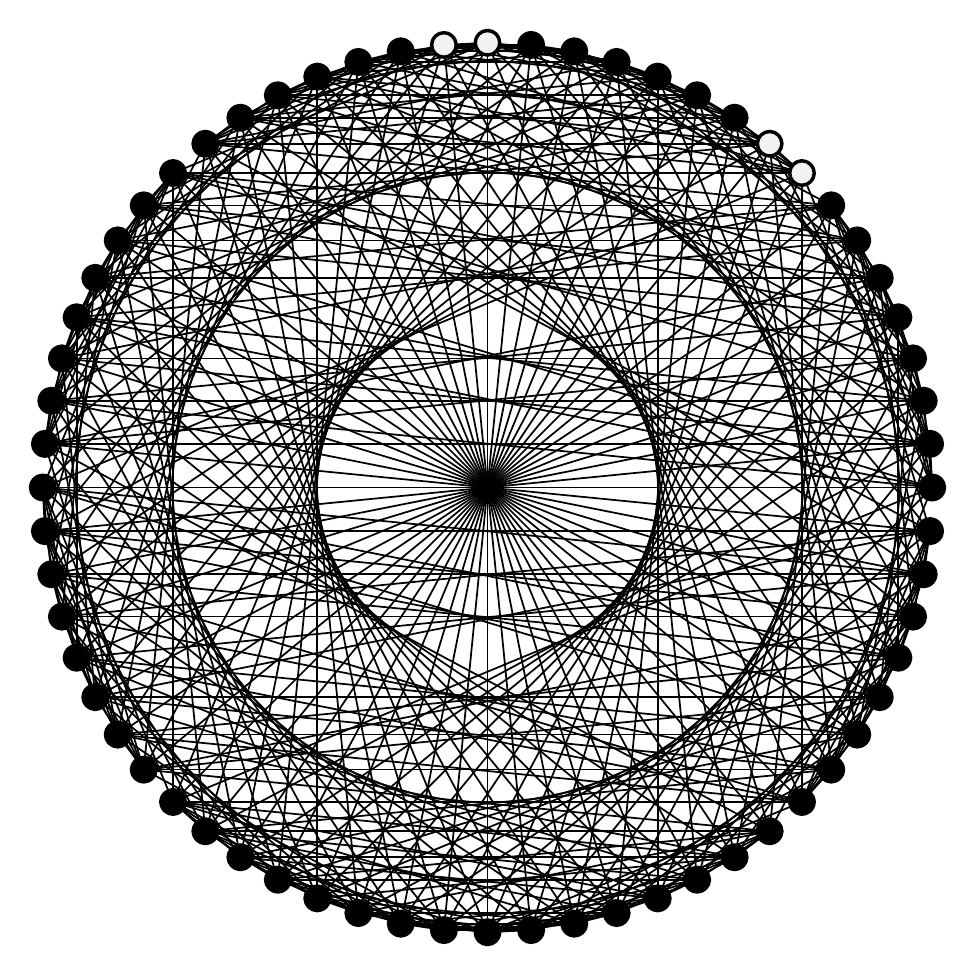}
  \caption{The graph $X$ for $d+1=6$}
  \label{fig:evendim6}
\end{figure}

For the even case, we still let $d\geq 5$ be odd, and we will consider a Cayley graph for $G:=\Z_2^{d+1}$. Let $\C$ be defined as before and let \(X=X(G,\C')\) where
\[\C':=\C\cup\{e_{d+1},\,\,e_1+\dots+e_d\}.\]
\begin{thm}\label{thm:exeven}
  Let $\C'$ be as described above for some odd $d\geq5$, define $X=X(\Z_2^{d+1},\C')$ and let $H\leq \Z_2^{d+1}$ be the maximal strongly cospectral subgroup in $X$. Then $|H|\geq4$.
\end{thm}
The graph is shown in Figure \ref{fig:evendim6} in dimension $6$.
The proof of this is longer, more technical and requires more cases than for the odd dimensions, and it will therefore be omitted from this paper. The idea is the following.

%\proof
\noindent\textsl{Idea of proof. }
  We see that the sum of $\C'$ is $e_{d+1}$ which is therefore strongly cospectral to $0$. Using similar methods as in the odd case, we can show that if $d+1\equiv 2\pmod4,$ then $e_1+e_2+e_3$ is strongly cospectral to zero and if $d+1\equiv 0\pmod4$, then $e_4+\cdots+e_d$ is strongly cospectral to zero.
\qed

\begin{lem}\label{lem:four}
  Let $X:=X(\Z_2^d,\C)$ be a cubelike graph, let $m$ be an odd number and denote by $C_m$ the cycle graph with vertices $\{0,\dots, m-1\}$. If the vertices $0$ and $g$ are strongly cospectral in $X$ then $(0,0)$ and $(g,0)$ are strongly cospectral in the Cartesian product $Y :=X\square C_m$, and this is a Cayley graph for the group $\Z_2^d\times \Z_m$.
\end{lem}
\proof
  Denote the vertices of $Y$ by $(g,h)$ where $g\in V(X) = \Z_2^d$ and $h\in V(C_m) = \Z_m$. Note that $V(Y) = \Z_2^d\times \Z_m=: G$. Define $\C_0:= \{(c,0): c\in \C\}$ and let $\C':=\C_0\cup \{(0,1),(0,-1)\}$. It is easy to verify that $Y = X(G,\C')$.

  If $\chi_1$ and $\chi_2$ are characters of $\Z_2^d$ and $\Z_m$, respectively, then the function $\chi:G\to \Cx^*$ defined by $(g,h)\mapsto \chi_1(g)\chi_2(h)$ is a character of $G$, and because $2^d$ and $m$ are coprime, distinct pairs $(\chi_1,\chi_2)$ give distinct characters $\chi$ of $G$, and so every character of $G$ can be decomposed uniquely in this way.

  Suppose $g\in \Z_2^d$ is strongly cospectral to $0$ in $X$. Let $\psi,\phi$ be characters of $G$ such that $\psi(\C') = \phi(\C')$.
  By the above, we can write $\phi = \phi_1\phi_2$ and $\psi=\psi_1\psi_2$ where $\phi_1,\psi_1$ are characters of $\Z_2^d$ and $\phi_2,\psi_2$ are characters of $\Z_m$. Then we have
  \begin{align*}
    \psi(\C') & = \sum_{c'\in\C'}\psi(c') = \sum_{c\in \C}\psi_1(c) + \psi_2(1)+\psi_2(-1) = \psi_1(\C)+\zeta_m^k+\zeta_m^{-k},
  \end{align*}
  for some $k< m$, where $\zeta_m$ is a primitive $m$-th root of unity. Similarly, there is some $\ell< m$ such that $\phi(\C') = \phi_1(\C')+\zeta_m^\ell+\zeta_m^{-\ell}$. Since $\psi(\C') = \phi(\C')$, this implies that
  \[(\zeta_m^k +\zeta_m^{-k}) - (\zeta_m^\ell +\zeta_m^{-\ell}) = \psi_1(\C)-\phi_1(\C)\in\Z,\]
  but since $m$ is odd, the only possibility is zero. Therefore, $\psi_1(\C) = \phi_1(\C)$ and so by Theorem \ref{thm:abel-evals}, $\psi_1(g) = \phi_1(g)$. This implies that $\psi(g,0) = \phi(g,0)$ which again implies that $(g,0)$ is strongly cospectral to $(0,0)$ in $Y$.
\qed

\begin{thm}\label{thm:fourall}
  If $n$ is a positive integer divisible by $32$ then there exists a Cayley graph on $n$ vertices with strongly cospectral sets of size four.
\end{thm}
\proof
  We can write $n=2^dm$ where $d\geq 5$ and $m$ is odd. By Theorems \ref{thm:exodd} and \ref{thm:exeven}, there is a cubelike graph on $2^d$ vertices with strongly cospectral sets of size four. Then, by Lemma \ref{lem:four}, we can build a Cayley graph of $\Z_2^d\times\Z_m$ preserving the strongly cospectral sets of size four.
\qed

\section{Further work}
There are several questions that we have not been able to answer.

Is the upper bound on the size of strongly cospectral sets in cubelike graphs tight in all dimensions?

Despite some effort, we did not find any cubelike graphs with strongly cospectral sets of size eight, and in fact we did have a more general question: Are there any vertex-transitive graphs with strongly cospectral sets of size more than four? This has in fact been answered very recently by Sin in \cite{sin2022strongly}, where he constructs cubelike graphs with arbitrarily large strongly cospectral sets (note that these sets are still far from reaching our upper bound in high dimensions).

The bound derived in Section \ref{sec:multnorm} holds for all normal Cayley graphs, but we have only been able to apply it to very restricted classes of graphs to get explicit bounds. Can we find a lower bound on the largest multiplicity of an eigenvalue in other normal Cayley graphs, or even Cayley graphs in general?

What about non-normal Cayley graphs? Can we find an upper bound on the size of strongly cospectral sets in the non-normal case?

\section*{Acknowledgments}
We would like to thank Gordon Royle for providing us with a complete list of cubelike graphs on $32$ and $64$ vertices. We thank Ada Chan, Harmony Zhan and Xiaohong Zhang for helpful discussions, in particular Ada Chan for bringing this topic to our attention. We would further like to thank Peter Sin for pointing out a small mistake in a proof in an earlier version of this paper.
Finally, we acknowledge support of NSERC (Canada), Grant No.\ RGPIN-9439.

\appendix
\section{Connection sets}
As was mentioned in Example \ref{ex:tight}, there are precisely twelve cubelike graphs (up to isomorphism) on $32$ vertices with strongly cospectral sets of size four. Below we give the connection sets of the ones that have degree at most $15$. The other six graphs are their complements.

All calculations were done in \texttt{Sage}. Here, we think of the group as the multiplicative abelian group generated by \texttt{f0, f1, f2, f3, f4}, where each \texttt{fi} has order two. The connection sets are listed below.
\begin{enumerate}
  \item \texttt{[f0, f1, f2, f3, f4, f0*f4, f1*f4, f2*f4, f3*f4, f0*f1*f2*f3]}
  \item \texttt{[f0, f1, f2, f3, f4, f0*f1, f2*f3, f2*f4, f3*f4, f0*f2*f3*f4,\\ f1*f2*f3*f4]}
  \item \texttt{[f0, f1, f2, f3, f4, f0*f4, f1*f4, f2*f3, f2*f4, f3*f4,\\ f0*f1*f2*f3, f0*f1*f2*f3*f4]}
  \item \texttt{[f0, f1, f2, f3, f4, f1*f4, f2*f3, f2*f4, f3*f4, f0*f1*f4,\\ f1*f2*f3, f0*f1*f2*f3, f0*f2*f3*f4]}
  \item \texttt{[f0, f1, f2, f3, f4, f0*f1, f0*f2, f0*f3, f1*f4, f2*f4, f3*f4,\\ f0*f1*f2*f4, f0*f1*f3*f4, f0*f2*f3*f4]}
  \item \texttt{[f0, f1, f2, f3, f4, f0*f1, f0*f2, f0*f3, f0*f4, f1*f2, f1*f3,\\ f1*f4, f2*f3, f2*f4, f0*f1*f2*f3*f4]}
\end{enumerate}
Table \ref{tab:spectra} shows the spectra of the corresponding graphs.
\begin{table}[h]
  \centering
  \setlength{\tabcolsep}{4pt} % Default value: 6pt
  \renewcommand{\arraystretch}{1.4}
\begin{tabular}{>{\stepcounter{rowcount}\therowcount)}l|c|l}
\multicolumn{1}{r}{}
& Degree & Spectrum \\
\hline
& 10 &\(\left\{-6^{(1)},\, -4^{(4)},\, -2^{(8)},\, 0^{(8)},\, 2^{(6)},\, 4^{(4)},\, 10^{(1)}\right\}\)\\
& 11 &\(\left\{-5^{(3)},\, -3^{(6)},\, -1^{(8)},\, 1^{(8)},\, 3^{(4)},\, 5^{(2)},\, 11^{(1)} \right\}\)\\
& 12 & \(\left\{-6^{(2)},\,-4^{(3)},\, -2^{(8)},\, 0^{(8)},\, 2^{(6)},\, 4^{(4)},\, 12^{(1)} \right\}\)\\
& 13 & \(\left\{-5^{(4)},\, -3^{(5)},\, -1^{(8)},\, 1^{(8)},\, 3^{(4)},\, 5^{(2)},\, 13^{(1)} \right\}\)\\
& 14 &\(\left\{-6^{(2)},\, -4^{(4)},\, -2^{(7)},\, 0^{(8)},\, 2^{(6)},\, 4^{(4)},\, 14^{(1)}\right\}\)\\
& 15 & \(\left\{-5^{(4)},\, -3^{(6)},\, -1^{(7)},\, 1^{(8)},\, 3^{(4)},\, 5^{(2)},\, 15^{(1)}\right\}\)
\end{tabular}
\caption{Spectra of the Cayley graphs with connection sets (1)-(6)}
\label{tab:spectra}
\end{table}

\section{Code}
We used the following \texttt{Sage} code to determine the maximal strongly cospectral subgroups in the  cubelike graphs. The code can be run on any abelian group.

\lstset{language=Python,%
    %basicstyle=\color{myred}
    basicstyle=\footnotesize,
    breaklines=true,%
    %morekeywords={for},
    keywordstyle=\color{myred},%
    morekeywords=[2]{1}, keywordstyle=[2]{\color{mygreen}},
    identifierstyle=\color{black},%
    stringstyle=\color{mylilas},
    commentstyle=\color{mygreen},%
    showstringspaces=false,%without this there will be a symbol in the places where there is a space
    numbers=left,%
    numberstyle={\tiny \color{black}},% size of the numbers
    numbersep=9pt, % this defines how far the numbers are from the text
    emph=[1]{StrCospVx,end,break,list,sum},emphstyle=[1]\color{myblue}, %some words to emphasise
    frame=single
    %emph=[2]{word1,word2}, emphstyle=[2]{style},
}
\vspace{0.2cm}

\begin{lstlisting}[language=Python]
  def StrCospVx(G,C):
      '''
      Input: An abelian group G, and a list C, subset of G.
      Output: A list of vertices that are strongly cospectral to 0 in the Cayley graph X(G,C).
      '''
      G1 = []
      for g in G.list():
          if g.order() == 2:
              G1.append(g)
      D = G.dual_group()
      # We create a dictionary where the keys are eigenvalues and the values are lists of the corresponding characters
      CharDict = {}
      for chi in D.list():
          chiC = sum([chi(x) for x in C])
          if chiC in CharDict.keys():
              CharDict[chiC].append(chi)
          else:
              CharDict[chiC] = [chi]
      Vx = [G.list()[0]]
      # We iterate over the elements of G of order 2
      CharLists = list(CharDict.values())
      for g in G1:
          i = 0
          while i>-1:
              # For a specified eigenvalue, make a list of the values of the corresponding eigenvectors in g
              Vals = [chi(g) for chi in CharLists[i]]
              # If they are not all the same, we stop
              if Vals.count(Vals[0]) != len(Vals):
                  i = -1
              # If we have not checked all the eigenvalues we keep going
              elif i < len(CharDict) - 1:
                  i = i+1
              # Otherwise, g is strongly cospectral to zero
              else:
                  Vx.append(g)
                  i = -1
      return Vx
\end{lstlisting}

\vspace{1.2 cm}

\bibliographystyle{plain}
%\bibliography{Bibliography.bib}

\end{document}